\documentclass[a4paper, 11pt]{article}

\usepackage[small]{titlesec}
\usepackage{tikz-cd}
\usepackage{booktabs}
\usepackage[T1]{fontenc}
\usepackage[utf8]{inputenc}
\usepackage[english]{babel}
\usepackage[a4paper,top=4cm,bottom=3cm,left=2.5cm,right=2.5cm]{geometry}
\usepackage{amsfonts}
\usepackage{mathtools}
\usepackage{amsmath}
\usepackage{xcolor}
\usepackage{array}
\usepackage{lmodern}
\usepackage{empheq}
\usepackage{amsthm}
\usepackage{amssymb}
\usepackage{booktabs}
\usepackage{caption}
\usepackage{siunitx}
\usepackage{float}
\usepackage{dsfont}
\usepackage{amssymb}
\usepackage{graphicx}
\usepackage{stmaryrd}

\usepackage{color}
\usepackage{pifont}%
\usepackage{enumerate}
\usepackage{enumitem}
\setlist[description]{leftmargin=\parindent,labelindent=\parindent}
\usepackage{listings} 
\usepackage{amsmath}
\usepackage{amscd}
\usepackage{braket}
\usepackage{tikz}
\usepackage{calc}
\usepackage[hidelinks]{hyperref}
\usepackage{resmes}
\usepackage{mathtools}
\usepackage{bbm}

\newcommand*{\Nearrow}{\rotatebox[origin=c]{45}{\(\Longrightarrow\)}}
\newcommand*{\Searrow}{\rotatebox[origin=c]{315}{\(\Longrightarrow\)}}

\newcommand{\R}{\mathbb{R}}
\newcommand{\sph}{\mathbb{S}}

\newcommand{\tra}{\mathrm{tr}}

\newcommand{\rhobar}{\overline{\rho}}

\newtheorem{theorem}{Theorem}[section]
\newtheorem*{theorem*}{Theorem}
\newtheorem{corollary}[theorem]{Corollary}

\theoremstyle{definition}

\usepackage{pgfplots}
\pgfplotsset{compat=1.17}

\newcommand{\di}{\, \mathrm{d}}

\newcommand{\La}{\mathcal{L}}

\usepackage{tikz-cd} 

\newcommand\blfootnote[1]{%
	\begingroup
	\renewcommand\thefootnote{}\footnote{#1}%
	\addtocounter{footnote}{-1}%
	\endgroup
}

\title{A strengthening of the dimensional Brunn--Minkowski conjecture implies the (B)-conjecture}
\author{Sotiris Armeniakos and Jacopo Ulivelli }
\date{}

\begin{document}

\maketitle

\begin{abstract}
In this note, we prove that if a sufficiently regular even log-concave measure satisfies a certain stronger form of the dimensional Brunn-Minkowski conjecture, then it also satisfies the (B)-conjecture. Furthermore, we show that hereditarily convex measures satisfy the aforementioned strengthened form, therefore providing an alternative proof of a recent result by Cordero-Erausquin and Eskenazis stating that a hereditarily convex measure satisfies both conjectures.
		
	\end{abstract}

\section{Introduction}

A measure $\mu$ on $\R^n$ is said to be log-concave if it satisfies the inequality
\begin{equation}\label{eq:log_concave_measure}
    \mu((1-t)K+tL)\geq \mu(K)^{1-t}\mu(L)^t
\end{equation}
for every convex and compact sets $K,L \subset \R^n$ and every $t \in [0,1]$, where $K+L$ is the Minkowski sum of $K$ and $L$. Among the central open problems in the study of log-concave measures are the dimensional Brunn--Minkowski conjecture and the (B)-conjecture.

The dimensional Brunn--Minkowski conjecture asks whether the dimensional concavity of the Lebesgue measure extends to all log-concave measures under symmetry assumptions. More precisely, if a log-concave measure $\mu$ is even, that is, $\mu(A)=\mu(-A)$ for every $\mu$-measurable set $A \subset \R^n$, and if $K,L$ are origin-symmetric convex sets, then the conjecture asserts that
\begin{equation}\label{eq:dim_BM}\tag{\rm dim-BM}
    \mu((1-t)K+tL)^{\frac{1}{n}}\geq (1-t)\mu(K)^{\frac{1}{n}}+t\mu(L)^{\frac{1}{n}}
\end{equation}
holds for every $t \in [0,1]$. Notice that \eqref{eq:dim_BM} in the case of the Lebesgue measure reduces to the celebrated Brunn--Minkowski inequality (see, e.g., \cite[Section~7]{SchneiderConvexBodiesBrunn2013}), which does not require any symmetry for $K$ and $L$. The dimensional Brunn--Minkowski conjecture was first stated by Gardner and Zvavitch \cite{gardner2010gaussian} for the Gaussian measure, a case settled by Eskenazis and Moschidis \cite{dim_gauss}, while Cordero-Erausquin and Rotem \cite{rot_B} proved \eqref{eq:dim_BM} for a broader class of measures, including sufficiently regular rotationally invariant log-concave measures. Nevertheless, the conjecture remains open for general even log-concave measures. See \cite{armen_uli, cordero2025concavity} for some recent developments and further references. 

The (B)-conjecture, seemingly independent of \eqref{eq:dim_BM}, concerns the behaviour of suitable dilates of convex bodies with respect to a log-concave measure. It asks whether, for an even log-concave measure $\mu$ on $\mathbb{R}^n$, the function
\begin{equation}\label{eq:B_conj} \tag{B}
    t \mapsto \mu(e^t K) \,, t \in \R, 
\end{equation}
is log-concave for every origin-symmetric convex $K \subset \R^n$. This problem was put forth by Banaszczyk and popularised by Lata{\l}a \cite{Latala2002ICM}. For the case of the Gaussian measure, this question was affirmatively answered by Cordero-Erausquin, Fradelizi, and Maurey \cite{Cordero-BGaussian}. A previous proof by Brascamp and Lieb \cite{BL_old}, which remained unnoticed for decades, recently resurfaced in \cite{cordero2025concavity}. Cordero-Erausquin and Rotem \cite{rot_B} later proved that the conjecture holds for the same general class of measures as \eqref{eq:dim_BM} above. More recently, Cordero-Erausquin and Eskenazis \cite{cordero2025concavity} established the conjecture for so-called hereditarily convex measures (see later \eqref{eq:hered_conv}), of which the measures treated in \cite{rot_B} are a subset. 

In \cite{cordero2025concavity} it is shown that hereditarily convex measures satisfy both \eqref{eq:dim_BM} and \eqref{eq:B_conj}, extending earlier work of Kolesnikov and Livshyts \cite{Kol_Liv_Gardner-Zvavitch}, who proved that, for sufficiently regular log-concave measures, hereditary convexity implies \eqref{eq:dim_BM}. Although no direct implication between \eqref{eq:dim_BM} and \eqref{eq:B_conj} is currently known, we show in this note that a suitable strengthening of \eqref{eq:dim_BM} yields \eqref{eq:B_conj}. Furthermore, we prove that this strengthening is implied by hereditary convexity, thereby providing an alternative route to some of the results in \cite{cordero2025concavity}. Our approach relies on a new method introduced in \cite{armen_uli}, combined with an application of the 
 $L_{2}$ method developed by Kolesnikov and Milman \cite{Kol_Mil_BL, Kol_Mil_P}, which was in turn inspired by \cite{Cordero-BGaussian}. 

Before presenting our results, we introduce some necessary notation. By a classical characterization due to Borell \cite{Borell}, a measure $\mu$ on $\R^n$ is log-concave if and only if there exists a convex function $u: \R^n \to \R \cup \{ \infty \}$ such that $\di \mu(x)=e^{-u(x)}\di x$. Throughout this paper, we will always consider $u \in C^2(\R^n)$ such that the Hessian matrix of $u$, denoted by $\nabla^2 u$, is positive definite. 

Let $K \subset \R^n$ be a compact and convex set with non-empty interior, such that its boundary $\partial K$ is a manifold of class $C^2$ with strictly positive principal curvatures. We shall say that such a set $K$ is of class $C^2_+$. Adopting the notation introduced by Livshyts \cite{GalynaStability}, the concavity power of $K$ with respect to $\mu$ is the maximal exponent $p(\mu,K)$ such that
\begin{equation}\label{eq:power_conc_gal}
    \left. \frac{\di^2}{\di t^2} \mu\bigl(K(\rho,t)\bigr)^{p(\mu,K)} \right|_{t=0}\leq 0\quad \text{for every } \rho \in \mathcal{A},
\end{equation}
where $K(\rho,t)$ denotes a perturbation of $K$ generated by a suitable function $\rho : \partial K \to \mathbb{R}$, and $\mathcal{A}$ is a class of admissible perturbations. We will elaborate more on this point in the next section. By a standard density argument, it can be shown that \eqref{eq:dim_BM} holds for $\mu$ if and only if $p(\mu,K) \geq 1/n$ for every origin-symmetric $K$ of class $C^2_+$, see \cite[Section~3.2]{GalynaStability} for a detailed exposition. Similarly, when studying the validity of \eqref{eq:B_conj}, we may, by a standard density argument, restrict our attention to compact sets with sufficiently regular boundaries.

Given a non-empty compact and convex set $K \subset \R^n$, we recall that its support function is defined as $h_K(x)=\sup_{y \in K} \langle x, y \rangle$. Moreover, if $K$ is of class $C^2_+$, we denote by $\nu_K: \partial K \to \sph^{n-1}$ its Gauss map. Throughout this note, given a compact and convex set $K$ with non-empty interior and a measure $\mu$ with continuous density with respect to the Lebesgue measure on $K$ (as will always be the case here), we ubiquitously use the symbol $\di \mu$ both for integration on $K$ and on $\partial K$. In the latter case, $\mu$ is understood as its restriction to the $(n-1)$-dimensional Hausdorff measure on $\partial K$. We are now ready to state our first main result.
\begin{theorem}\label{thm:main_1}
    Let $u \in C^2(\R^n)$ be an even function such that $\nabla^2 u$ is positive definite, and let $\mu$ be the associated log-concave measure with density $\di \mu(x)=e^{-u(x)}\di x$. Suppose that for every origin-symmetric set $K \subset \R^n$ of class $C^2_+$,  \begin{equation}\label{eq:strong_dim_BM}
        \frac{\mu(K)}{\int_{\partial K} h_K(\nu_K) \di \mu} \leq p(\mu,K). 
    \end{equation}
    Then, for every compact and convex and origin-symmetric set $K$ 
    \[ t \mapsto \mu(e^t K) \text{ is log-concave for } t \in \mathbb{R}. \]
    In other words, if $\mu$ satisfies \eqref{eq:strong_dim_BM} for every origin-symmetric $K$ of class $C^2_+$, $\mu$ satisfies the \eqref{eq:B_conj}-conjecture as well.
\end{theorem}\noindent
We say that the even log-concave measure $\mu$ satisfies the strong dimensional Brunn--Minkowski condition if \eqref{eq:strong_dim_BM} holds for every origin-symmetric set $K$ of class $C^2_+$. The fact that this condition implies the dimensional Brunn--Minkowski conjecture \eqref{eq:dim_BM} follows immediately from an integration by parts, together with the symmetry of the problem. See Corollary \ref{cor:strong_then_dim} later.

\medskip

To state our second result, we introduce some additional notation. Given $u$ and $\mu$ as in Theorem \ref{thm:main_1}, for every $\psi \in C^2(\R^n)$ we define the $\mu$-weighted Laplacian as $\La_\mu(\psi)=\Delta \psi-\langle \nabla u, \nabla \psi \rangle$. Moreover, we say that a measure $\eta$ is log-concave with respect to $\mu$ if $\eta$ is absolutely continuous with respect to $\mu$ and its corresponding density  $\di \eta / \di \mu$ is a log-concave function. The notion of hereditary convexity was formally introduced in \cite{ cordero2025concavity}. We restrict our treatment to the case of log-concave measures: A log-concave measure $\mu$ is hereditarily convex if it is even, and for every even measure $\eta$ which is log-concave with respect to $\mu$ and every even $\psi$ such that $\La_\mu(\psi) \in L_2(\eta)$, one has 
\begin{equation}\label{eq:hered_conv}
    \int_{\R^n} ||\nabla^2 \psi||^2_{\rm HS}+\langle \nabla^2 u \nabla \psi, \nabla \psi \rangle \di \eta \geq \frac{\left(\int_{\R^n} \La_\mu(\psi)\di \eta\right)^2}{\int_{\R^n} \La_\mu\left(\frac{|x|^2}{2}\right) \di \eta}.
\end{equation}
Here, $||\cdot||_{\rm HS}$ denotes the Hilbert--Schmidt norm of matrices, and $|x|^2=\langle x, x \rangle$. As previously mentioned, it is known that the class of hereditarily convex measures includes radially symmetric log-concave measures. It remains an open question of major interest whether every even log-concave measure is hereditarily convex; compare \cite[Question~16]{cordero2025concavity}. Our second result reads as follows.
\begin{theorem}\label{thm:main_2}
    Consider $\mu$ as in Theorem \ref{thm:main_1}. If $\mu$ is hereditarily convex, then $\mu$ satisfies the strong dimensional Brunn--Minkowski condition.
\end{theorem}\noindent
We shall prove, in fact, that the strong dimensional Brunn--Minkowski condition follows when $\mu$ is hereditarily convex only with respect to uniform measures on non-empty compact and convex subsets of $\R^n$. 

\medskip

To summarize this note, the results presented above show that if $\mu$ is an even-log concave measure as in Theorem \ref{thm:main_1}, the implications in \cite[(62)]{cordero2025concavity} can be refined as follows:
\begin{eqnarray}\notag
 & &    \textrm{$\mu$ verifies \eqref{eq:dim_BM}} \notag \\
 & \Nearrow &  \notag\\ \notag
 \textrm{$\mu$ is hereditarily convex} \Rightarrow \textrm{$\mu$ satisfies the strong dim-BM condition}  & &  \\
 & \Searrow  & \notag \\
 & &  \textrm{$\mu$ verifies \eqref{eq:B_conj}}. \notag
\end{eqnarray}

\section{Preliminaries and proofs of the results}

Before starting with the proofs, we introduce and recall some auxiliary notions and results. Consider a compact and convex set $K \subset \R^n$ of class $C^2_+$. First, we elaborate on the concavity power defined in \eqref{eq:power_conc_gal}. There, the perturbations $K(\rho,t)$ are defined as follows: As $K$ is of class $C^2_+$, it is well-known that for every $\rho \in C^2(\partial K)$ there exists $\varepsilon>0$ such that $h_t=h_K+t \rho(\nu^{-1}_K)$ are the support functions of compact and convex sets of class $C^2_+$ for every $|t|<\varepsilon$. We consider the sets $K(\rho,t)$ corresponding to the support functions $h_t$. Therefore, if $K$ is of class $C^2_+$, the family of admissible perturbations is, in fact, $C^2(\partial K)$. See also \cite[Section 3.2]{Kol_Mil_Lp_Brunn} and, for more details on the theory of convex sets, the monograph by Schneider \cite{SchneiderConvexBodiesBrunn2013}.

Given a set $K \subset \R^n$ of class $C^2_+$, we denote its second fundamental form by ${\rm II}$ and the tangential gradient on $\partial K$ by $\nabla_{\partial K}$. Consider, moreover, a function $u \in C^2(\R^n)$ with strictly positive Hessian matrix and the measure $\mu$ with density $\di \mu(x)=e^{-u(x)}\di x$. We write the weighted mean curvature of $\partial K$ as ${\rm H}_\mu=\tra({\rm II })- \langle \nabla u, \nu_K \rangle$. With a slight abuse of notation, we shall use $\langle \cdot, \cdot \rangle$ to denote also the restriction of the scalar product of $\R^n$ to the tangent spaces of $\partial K$. We write $H^1(\partial K, \mu)$ for the Sobolev space of functions on $\partial K$ with respect to the measure $\mu$. The following statement, proved in \cite[Theorem 1.3]{armen_uli}, provides an analytic representation for the concavity power $p(\mu,K)$.
\begin{theorem}\label{thm:minimal_power}
    Consider $\mu$ such that it satisfies the assumptions of Theorem~\ref{thm:main_1}, and a compact and convex origin-symmetric set $K \subset \R^n$ of class $C^2_+$. Then there exists a unique weak solution $\overline{\rho} \in H^{1}(\partial K,\mu)$ to the equation
\begin{equation}\label{eq:PDE_divergence_form_intro}
-\nabla_{\partial K} \cdot ({\rm II}^{-1}\nabla_{\partial K}\rho)+\langle \nabla_{\partial K}u,{\rm II}^{-1}\nabla_{\partial K}{\rho} \rangle -{\rm H}_\mu\,\rho+\frac{1}{\mu(K)}\int_{\partial K}\rho\,\di\mu=1\qquad\text{on }\partial K,
\end{equation}
and it satisfies 
\begin{equation}\label{eq:new_concavity expression_intro}
p(\mu,K) = \frac{\mu(K)}{\int_{\partial K} \overline{\rho}\,\di\mu}.
\end{equation}
Moreover, $\rhobar$ is even.
\end{theorem}\noindent
For later use, notice that \eqref{eq:new_concavity expression_intro} allows one to rewrite \eqref{eq:strong_dim_BM} as 
\begin{equation}\label{eq:equivalent_strong_BM}
    \int_{\partial K} h_K(\nu_K) \di \mu \geq \int_{\partial K} \rhobar \di \mu.
\end{equation}
Going back to condition \eqref{eq:strong_dim_BM}, observe that integration by parts yields the identity (see later \eqref{eq:parts_support} for a proof) \[ \int_{\partial K} h_K(\nu_K) \di \mu = n\mu(K)-\int_K \langle \nabla u, x \rangle\di \mu. \] If we plug this expression in \eqref{eq:strong_dim_BM}, we infer
\begin{equation}\label{eq:ezpz}
    p(\mu,K) \geq \frac{1}{n - \frac{1}{\mu(K)}\int_K \langle \nabla u, x \rangle \di \mu} \geq \frac{1}{n},
\end{equation}
where, for the last inequality, we used that whenever $u$ is even and convex, then $\int_K \langle \nabla u, x \rangle \di \mu \geq 0$. As anticipated in the introduction, one has the following immediate consequence.
\begin{corollary}\label{cor:strong_then_dim}
    Consider a log-concave measure $\mu$ as in Theorem \ref{thm:main_1}. If $\mu$ satisfies \eqref{eq:strong_dim_BM} for every compact, convex and origin-symmetric set of class $C^2_+$, then \eqref{eq:dim_BM} holds true.
\end{corollary}\noindent
As a further corollary of \eqref{eq:ezpz}, considering $tK$ for $t \to 0^+$ shows that the infimum of $p(\mu,K)$ on all the possible choices of $K$ cannot give, for $p(\mu,K)$, a lower bound larger than $1/n$. In fact, as mentioned, for example, in \cite{Kol_Liv_Gardner-Zvavitch}, an analogous limit argument shows that $1/n \geq \inf_{K} p(\mu,K)$.

\medskip

We now turn our attention to the PDE \eqref{eq:PDE_divergence_form_intro}. For $\rho \in C^2(\partial K)$ let us define the operator
\begin{equation}\label{eq:definition_operator}
    \mathcal{E}(\rho)=-\nabla_{\partial K} \cdot ({\rm II}^{-1}\nabla_{\partial K}\rho)+\langle \nabla_{\partial K}u,{\rm II}^{-1}\nabla_{\partial K}{\rho} \rangle -{\rm H}_\mu\,\rho+\frac{1}{\mu(K)}\int_{\partial K}\rho\,\di\mu.
\end{equation}
Then, $\mathcal{E}$ is a uniformly elliptic operator on $\partial K$ and $\rhobar$ is a weak solution of $\mathcal{E}(\rho)=1$.  We remark that if the support function of $K$ is assumed to be of class $C^3$ in $\R^n \setminus \{ 0 \}$, all the coefficients in \eqref{eq:definition_operator} are of class at least $C^1$. Under this assumption, the function $\rhobar$ is, in fact, a classical solution for \eqref{eq:PDE_divergence_form_intro}, as can be seen through classical elliptic regularity results. See, e.g., \cite{evans2022partial}. Now, observe that integration by parts yields the identity \[\int_{\partial K}\mathcal{E}(\rho)\rho \di \mu= \int_{\partial K} \langle {\rm II}^{-1} \nabla_{\partial K} \rho, \nabla_{\partial K} \rho \rangle \di \mu -\int_{\partial K} {\rm H}_\mu \rho^2 \di \mu+\frac{1}{\mu(K)}\left(\int_{\partial K} \rho \di \mu\right)^2.\] 
The right-hand side is non-negative thanks to the weighted Poincar\'e inequality (compare \cite{Kol_Mil_P}), which can be obtained by direct differentiation of \eqref{eq:log_concave_measure}. Therefore, 
\begin{equation}\label{eq:positive_def}
     \int_{\partial K}\mathcal{E}(\rho)\rho \di \mu \geq 0
\end{equation}
for every $\rho \in C^2(\partial K)$. We remark that this inequality holds for less regular functions. Additionally, we remark that the action of $\mathcal{E}$ on the support function, whenever $h_{K}$ is of class $C^2$ in $\R^n \setminus \{ 0 \}$, is
\begin{equation}\label{eq:action_support}
\mathcal{E}(h_{K}(\nu_{K})) = 1+\langle \nabla u, x \rangle -\frac{1}{\mu(K)}\int_{K} \langle \nabla u,x \rangle \di \mu.
\end{equation}
This equality can be found in \cite[(67)]{armen_uli}. It follows from direct computations, which we omit for the sake of brevity.

Next, we recall a local formulation of the (B)-conjecture. Differentiating twice $\mu(e^tK)$ for $\mu$ and $K$ as above (see, for example, \cite[(10)]{rot_B}, which in turn generalizes \cite[(4)]{Cordero-BGaussian}), shows that the log-concavity of the function $t \mapsto \mu(e^{t}K)$ at $t=0$ is equivalent to
\begin{equation}\label{eq:local_B}
    \frac{\int_K \langle\nabla u,x \rangle^2 \di \mu}{\mu(K)} - \left( \frac{\int_K \langle\nabla u,x \rangle \di \mu}{\mu(K)} \right)^2 \leq \frac{\int_K \langle (\nabla^2 u) x, x \rangle \di \mu}{\mu(K)} + \frac{\int_K \langle \nabla u, x \rangle \di \mu}{\mu(K)}.
\end{equation}
Thus, \eqref{eq:B_conj} is equivalent to the validity of \eqref{eq:local_B} for every compact and origin-symmetric convex set $K \subset \R^n$ with non-empty interior. 
For the $\mu$-weighted Laplacian $\La_{\mu}$, we recall the following useful integration by parts identity, which holds for any two functions $\varphi,\psi \in C^2(K)$:
\begin{equation}\label{eq:integration_by_parts}
\int_K \varphi \, \La_\mu (\psi) \di\mu = - \int_K \langle \nabla \varphi, \nabla \psi \rangle \di\mu+\int_{\partial K}\varphi \langle \nabla \psi, \nu_K \rangle \di\mu .
\end{equation}

Finally, we will make use of the following Reilly formula, proved in \cite[(2.3)]{Kol_Mil_BL} in the more general context of weighted Riemannian manifolds (compare also \cite[(21)]{cordero2025concavity}). Consider a function $\psi \in C^2(\R^n)$ and its normal derivative $\rho = \langle \nabla \psi, \nu_K \rangle$. Then
\begin{align}\label{eq:Reilly_formula}
\notag \int_{K} (\La_{\mu}  \psi)^2 \di \mu &= \int_{K} \| \nabla^2 \psi \|_{{\rm HS}}^2+\langle (\nabla^2 u) \nabla \psi, \nabla \psi \rangle \di \mu \\
&+ \int_{\partial K} {\rm H}_{\mu}\rho^2 + \langle {\rm II} \nabla_{\partial K} \psi, \nabla_{\partial K}\psi\rangle -2\langle \nabla _{\partial K}\psi, \nabla_{\partial K}\rho \rangle \di \mu. 
\end{align}

\medskip

We are now ready to prove our results.

\paragraph{Proof of Theorem \ref{thm:main_1}.}  We begin by following the argument in \cite[Proposition~5.3]{armen_uli}. Let $K \subset \R^n$ be a compact and convex origin-symmetric set and let $\rhobar$ be given by \eqref{eq:PDE_divergence_form_intro}. We may assume by approximation that, in addition, the support function of $K$ is of class $C^3$ with strictly positive principal curvatures and, thus, $\rhobar \in C^{2}(\partial K)$. Moreover, denote by $\eta$ the measure with density $\di \eta = \frac{1}{\mu(K)}\mathbbm{1}_K\di \mu$, where $\mathbbm{1}_K$ is the characteristic function of $K$ such that $\mathbbm{1}_K(x)=1$ if $x \in K$ and $0$ otherwise. Suppose that \eqref{eq:strong_dim_BM} holds, we shall show that \eqref{eq:local_B} is also satisfied. From \eqref{eq:action_support} and since $\mathcal{E}(\rhobar)=1$ we infer
\begin{align*}
\int_{\partial K} h_{K}(\nu_{K})-\rhobar\di \eta &= \int_{\partial K}(h_{K}(\nu_{K})-\rhobar)\mathcal{E}(\rhobar)\di \eta = \int_{\partial K}\mathcal{E}(h_{K}(\nu_{K})-\rhobar)\rhobar\di \eta \\
&=\int_{\partial K} \rhobar \langle \nabla u,x \rangle \di \eta - \int_{\partial K} \rhobar \di \eta \int_{K} \langle \nabla u,x \rangle \di \eta.
\end{align*}
Then, writing $f(x)=\rhobar(x) -h_{K}(\nu_{K}(x))$,
\begin{align}\label{eq:2.1.1}
         \notag \int_{\partial K} h_{K}(\nu_{K})-\rhobar\di \eta &= \underbrace{\int_{\partial K} f \langle \nabla u,x \rangle \di \eta - \int_{\partial K} f \di \eta \int_{K} \langle \nabla u,x \rangle \di \eta}_{({\rm I})} \\
         &+\underbrace{\int_{\partial K} h_{K}(\nu_{K}) \langle \nabla u,x \rangle \di \eta - \int_{\partial K} h_{K}(\nu_{K}) \di \eta \int_{K} \langle \nabla u,x \rangle \di \eta}_{({\rm I\! I})}.
\end{align}
We begin by computing the term $({\rm I})$. Notice that by \eqref{eq:action_support} and \eqref{eq:PDE_divergence_form_intro} we have
\[\mathcal{E}(f)=\mathcal{E}(\rhobar-h_{K}(\nu_K)) = - \langle \nabla u,x\rangle +\int_{K} \langle \nabla u,x \rangle \di \eta .\]
Thus,
\[({\rm I}) = -\int_{\partial K}f\mathcal{E}(f) \di \eta \leq 0 \,,\] 
where the last inequality is \eqref{eq:positive_def}. We proceed with the term $({\rm I\! I})$. Notice that since for every $x \in \partial K$ we have that $h_{K}(\nu_{K}(x))= \langle x,\nu_{K}(x)\rangle = \langle \nabla (|x|^2/2), \nu_K(x) \rangle$, an application of \eqref{eq:integration_by_parts} implies
\begin{equation}\label{eq:parts_support}
    \int_{\partial K}h_{K}(\nu_K)\di \eta = \int_{K}\La_{\mu}\left(\frac{|x|^2}{2}\right)\di \eta = n -\int_{K}\langle \nabla u,x \rangle \di \eta.
\end{equation}
Moreover, integrating again by parts shows
\begin{align*}
\int_{\partial K}\langle \nabla u,x \rangle h_{K}(\nu_K) \di \eta &= \int_{K}\La_{\mu}\left(\frac{|x|^2}{2}\right)\langle \nabla u, x\rangle \di \eta +\int_{K}  \langle \nabla \langle \nabla u,x\rangle , x\big \rangle \di \eta \\ 
&=n\int_{K}\langle \nabla u,x \rangle \di \eta -\int_{K}\langle \nabla u,x \rangle^2 \di \eta+\int_{K}\langle (\nabla^2u)x,x\rangle \di \eta+ \int_{K}\langle \nabla u,x \rangle \di \eta.
\end{align*}
Therefore, we infer
\begin{equation}\label{eq:eq:2.1.2}
    ({\rm I\! I}) = \int_{K}\langle (\nabla^2u)x,x\rangle \di \eta+\int_{K}\langle \nabla u,x\rangle\di \eta - \int_{K}\langle \nabla u,x \rangle^2\di \eta+\left(\int_{K} \langle \nabla u, x \rangle \di \eta\right)^2.
\end{equation}
If the strong dimensional Brunn--Minkowski condition holds, by \eqref{eq:equivalent_strong_BM} we have that $({\rm I})+({\rm I\! I}) \geq 0$ and, since $({\rm I})\leq 0$, we conclude that $({\rm I\! I}) \geq 0$. That is,
\[\int_{K}\langle (\nabla^2u)x,x\rangle \di \eta+\int_{K}\langle \nabla u,x\rangle\di \eta \geq \int_{K}\langle \nabla u,x \rangle^2\di \eta-\left(\int_{K} \langle \nabla u, x \rangle \di \eta\right)^2,\]
which is precisely \eqref{eq:local_B}. Since $K$ was arbitrary we infer that \eqref{eq:B_conj} holds, concluding the proof. 

\paragraph{Proof of Theorem \ref{thm:main_2}.} For this proof, we apply H\"ormander's $L_2$ method, following the approach introduced in \cite{Kol_Mil_BL, Kol_Mil_P}, which can be consulted for further details. See also the discussion in \cite[Section 2]{cordero2025concavity}. Given $K \subset \R^n$ of class $C^2_+$, we shall prove that \eqref{eq:strong_dim_BM} is satisfied. Consider the function $\rhobar$ solving \eqref{eq:PDE_divergence_form_intro} for the previous choice of $K$, and let $\psi$ be the unique solution to the Neumann problem 
\begin{equation}\label{eq:2.2.1}
\begin{cases}
\La_\mu \psi= \frac{1}{\mu(K)}\int_{\partial K} \overline{\rho}\di \mu &\text{in }K,\\[2pt]
\displaystyle \langle \nabla \psi, \nu_K \rangle=\overline\rho&\text{on }\partial K.
\end{cases}
\end{equation}
Such solution exists provided that the integral (with respect to $\mu$) of the interior data on $K$ is the same as the integral of the boundary data on $\partial K$, which is indeed the case here. By approximation, we can assume $h_{K}$ to be of class $C^3$ in $\R^n \setminus \{0 \}$, so that $\rhobar \in C^2(\partial K)$. Compare the discussion after \eqref{eq:definition_operator}. Thus, standard results in elliptic regularity theory (see, for example, \cite[Theorem 2.5]{Kol_Mil_BL}) yield that $\psi \in C^2(K)$. Consider the normalized measure $\eta$ as in the previous proof. The Cauchy--Schwarz inequality for positive quadratic forms implies that
\[ 0 \leq \langle {\rm II} \nabla_{\partial K} \psi, \nabla_{\partial K} \psi \rangle - 2 \langle \nabla_{\partial_K} \psi, \nabla_{\partial K} \rhobar \rangle+\langle {\rm II}^{-1} \nabla_{\partial K} \rhobar, \nabla_{\partial K} \rhobar \rangle.\]
Applying this to the Reilly formula \eqref{eq:Reilly_formula}, for $\mu$, $\psi$ and dividing by $\mu(K)$ implies
\begin{align}\label{eq:2.2.2}
\int_{\partial K} \langle {\rm II}^{-1} \nabla_{\partial K} \rhobar, \nabla_{\partial K} \rhobar \rangle \di \eta-\int_{\partial K}{\rm H}_{\mu}\rhobar^2\di \eta +\left(\int_{\partial K} \rhobar \di \eta\right)^2 \geq \int_{K} \|\nabla \psi\|^2_{{\rm HS}} + \langle (\nabla u)^2 \nabla \psi , \nabla \psi \rangle\di \eta,
\end{align}
Furthermore, since by \eqref{eq:PDE_divergence_form_intro} $\rhobar$ satisfies $\mathcal{E}(\rhobar) = 1$, multiplying with $\rhobar$ and integrating yields
\begin{equation}\label{eq:2.2.3}
\int_{\partial K} \langle {\rm II}^{-1} \nabla_{\partial K} \rhobar, \nabla_{\partial K} \rhobar \rangle \di \eta- \int_{\partial K}{\rm H}_{\mu}\rho^2\di \eta +\left(\int_{\partial K} \rhobar \di \eta\right)^2 = \int_{\partial K} \rhobar \di \eta.
\end{equation}
Notice that this implies that $\int_{\partial K} \rhobar \di \eta \geq 0$. Furthermore, by \cite[Corollary~3.4]{armen_uli} equality in the corresponding Poincar\'e inequality holds if and only if the function is identically zero. Consequently, we infer that $\int_{\partial K} \rhobar\di  \eta >0$, as the alternative $\rhobar \equiv 0$ would contradict \eqref{eq:PDE_divergence_form_intro}.
Now, notice that 
\[\int_{K}\mathcal{L}_{\mu}\psi \di \eta= \int_{K}\frac{1}{\mu(K)}\int_{\partial K} \rhobar \di \mu \di\eta = \int_{\partial K} \rhobar \di \eta,\]
by definition of the measure $\eta$. Also, by Theorem~\ref{thm:minimal_power}, $\rhobar$ is even and thus so is $\psi$. Furthermore, since $\int_{K}\La_{\mu}(\frac{x^2}{2})\di \eta = \int_{\partial K}h_{K}(\nu_{K})\di \eta$ by \eqref{eq:parts_support}, we conclude by the hereditary convexity assumption \eqref{eq:hered_conv} that
\[\int_{K} \|\nabla \psi\|^2_{{\rm HS}} + \langle (\nabla u)^2 \nabla \psi , \nabla \psi \rangle \di\nu\geq \frac{\left(\int_{\partial K} \rhobar \di \eta\right)^2}{\int_{\partial K} h_{K}(\nu_{K})\di \eta},\]
which, in combination with \eqref{eq:2.2.2}, \eqref{eq:2.2.3}, and the fact that $\int_{\partial K}\rhobar \di \eta >0$ implies
\[\int_{\partial K} \rhobar \di \eta \leq \int_{\partial K} h_{K}(\nu_{K})\di \eta.\]
Thus, \eqref{eq:strong_dim_BM} holds for an arbitrary compact, convex and origin-symmetric set $K$ and the strong dimensional Brunn--Minkowski condition holds for the measure $\mu$, concluding the proof. 

\medskip

\paragraph{Acknowledgments.} The first-named author was supported by the Austrian Science Fund (FWF): 10.55776/P36344N. The second-named author was supported, in part, by the Austrian Science Fund (FWF): 10.55776/PAT3787224.

\footnotesize
\bibliography{references}
\bibliographystyle{siam}

\parbox[t]{8.5cm}{
Sotiris Armeniakos\\
Institut f\"ur Stochastik und Wirtschaftsmathematik\\
TU Wien\\
Wiedner Hauptstra{\ss}e 8-10/1046\\
1040 Wien, Austria\\
e-mail: sotirios.armeniakos@tuwien.ac.at}

\medskip

\parbox[t]{8.5cm}{
Jacopo Ulivelli\\
Institut f\"ur Diskrete Mathematik und Geometrie\\
TU Wien\\
Wiedner Hauptstra{\ss}e 8-10/1046\\
1040 Wien, Austria\\
e-mail: jacopo.ulivelli@tuwien.ac.at}
\end{document}